\documentclass[a4paper]{article}

\usepackage[latin1]{inputenc}
\usepackage[T1]{fontenc}
\usepackage[francais,english]{babel}

\usepackage{amsmath,amssymb,amsfonts,amsthm}
\usepackage{amsfonts}
\usepackage{amssymb}
\usepackage{latexsym}
\usepackage{comment}

\newtheorem{thm}{Theorem}
\newtheorem{Prop}{Proposition}
\newtheorem{Def}{Definition}
\newtheorem{Lem}{Lemma}

\author{
K. \textsc{Beauchard} 
\footnote{CMLA, ENS Cachan, CNRS, Universud, 
61 avenue du Pr\'esident Wilson, F-94230 Cachan, France,
email: Karine.Beauchard@cmla.ens-cachan.fr},
P. \textsc{Cannarsa}
\footnote{Università di Roma Tor Vergata, via della Ricerca Scientifica 1, 00133, Roma, Italy, 
email: cannarsa@axp.mat.uniroma2.it (corresponding author)},
R. \textsc{Guglielmi}
\footnote{Università di Roma Tor Vergata, via della Ricerca Scientifica 1, 00133, Roma, Italy, 
email: guglielm@axp.mat.uniroma2.it}
\thanks{This research has been performed in the framework of the GDRE CONEDP. The authors wish to thank Institut Henri Poincar\'e (Paris, France) for providing a very stimulating environment during the "Control of Partial and Differential Equations and Applications" program in the Fall 2010. This work was started during the C.I.M.E. course `Control of partial differential equations' (Cetraro, July 19--23, 2010).  The first author was partially supported by the ``Agence Nationale de la Recherche'' (ANR),
Projet Blanc C-QUID number BLAN-3-139579.}
}

\title{Null controllability of Grushin-type operators\\ in dimension two}
\date{}

\begin{document}

\maketitle

\begin{abstract}
We study the null controllability of the parabolic equation associated with the Grushin-type 
operator $A=\partial_x^2+|x|^{2\gamma}\partial_y^2\,, (\gamma>0),$ in the rectangle $\Omega=(-1,1)\times(0,1)$, under an additive control supported in the strip $\omega=(a,b)\times(0,1)\,, (0<a,b<1)$. 
We prove that the equation is null controllable in any positive time for $\gamma<1$, and that it fails to
be so for $\gamma>1$. In the transition regime $\gamma=1$, we show that both behaviors live together: 
a positive minimal time is required for null controllability.
Our approach is based on the fact that, thanks to the particular geometric configuration, null controllability is equivalent to the observability of the Fourier components of the solution of the adjoint system uniformly with respect to the frequency.
\end{abstract}



\section{Introduction}

\subsection{Main result}

We consider the Grushin-type equation
\begin{equation} \label{Grushin_eq}
\left\lbrace \begin{array}{ll}
\partial_t f - \partial_x^2 f - |x|^{2\gamma} \partial_y^2 f= u(t,x,y)1_{\omega}(x,y)& (t,x,y) \in (0, \infty)\times\Omega\,,\\
f(t,x,y)=0 & (t,x,y) \in (0, \infty)\times\partial \Omega\,,
\end{array}\right.
\end{equation}
where  $\Omega:=(-1,1) \times (0,1)$, $\omega \subset \Omega$, and $\gamma >0$.
Problem (\ref{Grushin_eq}) is a linear control system in which
\begin{itemize}
\item the state is $f$,
\item the control $u$ is supported in the subset $\omega$.
\end{itemize}
It is a degenerate parabolic equation, since the coefficient of $\partial_y^2f$ vanishes on the line $\{x=0\}$. We will investigate the null controllability of (\ref{Grushin_eq}).

\begin{Def}[Null controllability] Let $T>0$. System (\ref{Grushin_eq}) is {\em null controllable in time} $T$ if,
for every $f_0 \in L^2(\Omega)$, there exists $u \in L^{2}((0,T)\times \Omega ;\mathbb{R})$ such that the solution of
\begin{equation} \label{Grushin_CYpb}
\left\lbrace \begin{array}{ll}
\partial_t f - \partial_x^2 f - |x|^{2\gamma} \partial_y^2 f= u(t,x,y)1_{\omega}(x,y) &(t,x,y) \in (0,T)\times\Omega\,,\\
f(t,x,y)=0 & (t,x,y) \in (0,T)\times \partial \Omega\,,\\
f(0,x,y)=f_0(x,y) & (x,y) \in \Omega\,,
\end{array}\right.
\end{equation}
satisfies $f(T,\cdot,\cdot)=0$.

System (\ref{Grushin_eq}) is {\em null controllable} if there exists $T>0$ such that it is
null controllable in time $T$.
\end{Def}

The main result of this paper is the following one.

\begin{thm} \label{Main-thm}
Let $\omega=(a,b) \times (0,1)$, where $0<a<b<1$.
\begin{enumerate}
\item If $\gamma \in (0,1)$, then system (\ref{Grushin_eq}) is null controllable in any time $T>0$.
\item If $\gamma=1$, then there exists $T^*>0$ such that
\begin{itemize}
\item for every $T>T^*$ system (\ref{Grushin_eq}) is null controllable in time $T$,
\item for every $T<T^*$ system (\ref{Grushin_eq}) is not null controllable in time $T$.
\end{itemize}
\item If $\gamma >1$, then (\ref{Grushin_eq}) is not null controllable.
\end{enumerate}
\end{thm}

By duality, the null controllability of (\ref{Grushin_eq}) is equivalent to an observability inequality for the adjoint system 
\begin{equation} \label{adjointL2_ssCI}
\left\lbrace \begin{array}{ll}
\partial_t g - \partial_x^2 g - |x|^{2\gamma} \partial_y^2 g=0 & (t,x,y) \in (0,\infty)\times\Omega\,,\\
g(t,x,y)=0 & (t,x,y) \in (0, \infty)\times\partial \Omega\,.
\end{array}\right.
\end{equation}

\begin{Def} 
[Observability] Let $T>0$. System (\ref{adjointL2_ssCI}) is {\em observable in $\omega$ in time} $T$ if
there exists $C>0$ such that,
for every $g_{0} \in L^{2}(\Omega)$, the solution of
\begin{equation} \label{adjointL2}
\left\lbrace \begin{array}{ll}
\partial_t g - \partial_x^2 g - |x|^{2\gamma} \partial_y^2 g=0 & (t,x,y) \in (0,T)\times\Omega\, ,\\
g(t,x,y)=0 & (t,x,y) \in (0,T)\times \partial \Omega\, ,\\
g(0,x,y)=g_0(x,y) & (x,y) \in \Omega\, ,
\end{array}\right.
\end{equation}
satisfies
$$\int_{\Omega} |g(T,x,y)|^{2} dx dy 
\leqslant C \int_{0}^{T} \int_{\omega} |g(t,x,y)|^{2} dx dy dt\, .$$

System (\ref{adjointL2_ssCI}) is {\em observable in $\omega$} if there exists $T>0$ such that it is observable in $\omega$ in time $T$.
\end{Def}

\begin{thm} \label{thm:OI}
Let $\omega=(a,b) \times (0,1)$, where $0<a<b<1$.
\begin{enumerate}
\item If $\gamma \in (0,1)$, then system (\ref{adjointL2}) is observable in $\omega$ in any time $T>0$.
\item If $\gamma=1$, then there exists $T^*>0$ such that
\begin{itemize}
\item for every $T>T^*$  system (\ref{adjointL2}) is observable in $\omega$ in time $T$,
\item for every $T<T^*$  system (\ref{adjointL2}) is not observable in $\omega$ in time $T$.
\end{itemize}
\item If $\gamma >1$, then system (\ref{adjointL2}) is not observable in $\omega$.
\end{enumerate}
\end{thm}

\subsection{Motivation and bibliographical comments}

\subsubsection{Null controllability of the heat equation}

The null and approximate controllability 
of the heat equation are essentially well understood subjects
for both linear and semilinear equations, and
for bounded or unbounded domains
(see, for instance,
\cite{Doubova-EFC-MGB-EZ}, 
\cite{Fabre-Puel-Zuazua},
\cite{Cara-Zuazua-ChaleurNL},
\cite{Cara-Zuazua-ChaleurLineaire}, 
\cite{Cara-Zuazua-ChaleurBV},
\cite{Burgos-Teresa},
\cite{Imanuvilov-Yamamoto},
\cite{Lebeau-Robbiano}, 
\cite{Lopez-Zuazua}, 
\cite{Miller-unbounded},
\cite{Miller},
\cite{Zuazau-chaleurNLapp},
\cite{Zuazau-chaleurdf}).
Let us summarize one of the existing main results. Consider the linear heat equation
\begin{equation} \label{Chaleur-lineaire}
\left\lbrace\begin{array}{ll}
\partial_t f -\Delta f = u(t,x)1_{\omega}(x) &  (t, x) \in (0,T)\times\Omega\,,
\\
f=0 &\text{  on } (0,T) \times \partial \Omega\,,
\\
f(0,x)=f_{0}(x) & x\in \Omega\,,
\end{array} \right.
\end{equation}
where $\Omega$ is an open subset of $\mathbb{R}^{d}$, $d \in \mathbb{N}^{*}$, and
$\omega$ is a subset of $\Omega$.
The following theorem is due, for the case $d =1$, to H. Fattorini and D. Russell 
\cite[Theorem~3.3]{Fattorini-Russel}, and, for $d \geqslant 2$,
to O. Imanuvilov 
\cite{Imanuvilov-parabolic-boundary}, \cite{Imanuvilov-parabolic} 
(see also the book \cite{Fursikov-Imanuvilov-186} by A. Fursikov and O.Imanuvilov)
and  G. Lebeau and L. Robbiano \cite{Lebeau-Robbiano}.

\begin{thm} \label{Heat-borne}
Let us assume that $\Omega$ is bounded, of class $C^{2}$ and connected, $T>0$,
and $\omega$ is a nonempty open subset of $\Omega$. 
Then the control system (\ref{Chaleur-lineaire}) is null controllable in time $T$.
\end{thm}
So, the heat equation on a smooth bounded domain is null controllable
\begin{itemize}
\item in arbitrarily small time;
\item with an arbitrarily small control support $\omega$.
\end{itemize}

It is natural to ask whether null controllability also holds for degenerate parabolic equations such as (\ref{Grushin_eq}). 
Let us compare the known results for the heat equation with the results proved in this article.
The first difference concerns the geometry of $\Omega$ and $\omega$: a more restrictive configuration is assumed in   
Theorem~\ref{Main-thm} than in Theorem~\ref{Heat-borne}.
The second difference concerns the structure of the controllability results.
Indeed, while the heat equation is null controllable in arbitrarily small time, the same result holds for the Grushin equation only when degeneracy is not too strong (i.e. $\gamma \in (0,1)$). On the contrary, when degeneracy is too strong (i.e. $\gamma>1$), null controllability does not hold any more.
Of special interest is the transition regime ($\gamma=1$), where the `classical' Grushin operator appears: here, both behaviors live together, and a positive minimal time is required for the null controllability.

\subsubsection{Boundary-degenerate parabolic equations}

The null controllability of parabolic equations degenerating on the boundary of the domain in one space dimension 
is well-understood, much less so in higher dimension.  Given $0<a<b<1$ and $\alpha>0$, let us consider the 1D equation
$$\partial_t w + \partial_x (x^\alpha \partial_x w) = u(t,x) 1_{(a,b)}(x)\,,\qquad (t,x) \in (0,\infty)\times(0,1)\,,$$
with suitable boundary conditions. Then, it can be proved that null controllability holds
if and only if $\alpha \in (0,2)$ (see \cite{Cannarsa-V-M-ADE,Cannarsa-V-M-SIAM}), while, for $\alpha\ge 2$, the best result one 
can show is \textquotedblleft regional null controllability\textquotedblright  (see \cite{CMV-regional}), which consists in controlling the solution within the domain of influence of the control. Several extensions of the above results are available in one space dimension, see
 \cite{Ala-Can-Fra, Martinez-Vancost-JEE-2006} for equations in divergence form, \cite{Can-Fra-Roc_1,Can-Fra-Roc_2} for nondivergence form operators, and \cite{Can-deT, Flo-deTer} for cascade systems. Fewer results are available for multidimensional problems, mainly in the case of two dimensional parabolic operators which simply degenerate in the normal direction to the boundary of the space domain, see \cite{Cannarsa-V-M-CRAS}. As in the above references, also for the Grushin equation null controllability holds if and only if the
degeneracy is not too strong ($\gamma \in (0,1]$).

\subsubsection{Parabolic equations degenerating inside the domain}

In \cite{Martinez-Vancostenoble-Raymond}, the authors study linearized Crocco type equations
$$\left\lbrace\begin{array}{ll}
\partial_{t}f + \partial_{x}f - \partial_{vv}f=u(t,x,v) 1_{\omega}(x,v) &
(t,x,v) \in (0,T) \times (0,L) \times (0,1)\,,
\\
f(t,x,0)=f(t,x,1)=0 &
(t,x) \in (0,T) \times (0,L)\,,
\\
f(t,0,v)=f(t,L,v)
&
(t,v) \in (0,T)  \times (0,1)\,.
\end{array}\right.$$
For a given open subset $\omega$ of $ (0,L) \times (0,1)$,
they prove regional null controllability. 
Notice that, in the above equation, diffusion (in $v$) and transport (in $x$) are decoupled.

In \cite{Kolm}, the authors study the Kolmogorov equation
\begin{equation} \label{Kolm}
\partial_{t}f +v \partial_{x}f - \partial_{vv}f=u(t,x,v)1_\omega(x,v)\,, \quad (x,v) \in (0,1)^2\,,
\end{equation}
with periodic type boundary conditions.
They prove null controllability in arbitrarily small time,
when the control region $\omega$ is a strip, parallel to the $x$-axis.
We note that the above Kolmogorov equation degenerates  on the whole space domain, unlike Grushin's equation.
However, differently from the linearized Crocco equation, transport (in $x$ at speed $v$) and diffusion (in $v$) are coupled.
This is why the null controllability results are also different for these equations.

\subsubsection{Unique continuation and approximate controllability}
In this paper we will not directly address approximate controllability, which is another interesting problem in control theory.
It is well-known that, for evolution equations, approximate controllability can be equivalently formulated as 
unique continuation (see \cite{Zab}). The unique continuation problem for the Grushin-type operator 
$$A=\partial_x^2+|x|^{2\gamma}\partial_y^2$$
has been widely investigated. In particular,  in \cite{Garofalo2} (see also the references therein) 
unique continuation is proved  for every $\gamma>0$ and every open set $\omega$.

\subsubsection{Null controllability and hypoellipticity}

It could be interesting to analyze the connections between null controllability and hypoellipticity. We recall that 
a linear differential operator $P$ with $C^\infty$ coefficients in an open set $\Omega  \subset \mathbb{R}^n$
is called hypoelliptic if, for every distribution $u$ in $\Omega$, we have
$$\text{sing} \text{  } \text{supp} u = \text{sing} \text{  } \text{supp} Pu\, ,$$
that is, $u$ must be a $C^\infty$ function in every open set where so is $Pu$.
The following sufficient condition (which is also essentially necessary)
for hypoellipticity is due to H\"ormander (see \cite{HormHypoel}).
\begin{thm}
Let $P$ be a second order differential operator of the form
$$P=\sum_{j=1}^r X_j^2 + X_0 + c\,,$$
where $X_0,...,X_r$ denote first order homogeneous differential operators in an open set $\Omega \subset \mathbb{R}^n$
with $C^\infty$ coefficients, and $c \in C^\infty(\Omega)$.
Assume that there exists $n$ operators among 
$$X_{j_1}\, ,\ [X_{j_1},X_{j_2}]\,,\ [X_{j_1},[X_{j_2},X_{j_3}]]\, ,\dots ,\
[X_{j_1},[X_{j_2},[X_{j_3},[...,X_{j_k}]...]]]\,,$$ 
where $j_i \in \{0,1,...,r\}$, which are linearly independent at any given point in $\Omega$. Then, $P$ is hypoelliptic.
\end{thm}

H\"ormander's condition is satisfied by the Grushin operator $A=\partial_x^2+|x|^{2\gamma}\partial_y^2$ for every $\gamma \in \mathbb{N}^*$ (for other values of $\gamma$, the coefficients are not $C^\infty$).
Indeed, set
$$X_1(x,y):=\left( \begin{array}{c}
1 \\ 0
\end{array} \right)\,,\quad
X_2(x,y):=\left( \begin{array}{c}
0 \\ x^\gamma
\end{array} \right)\,.$$
Then,
$$
[X_1,X_2](x,y)=\left( \begin{array}{c}
0 \\ \gamma x^{\gamma-1}
\end{array} \right)\, ,\ 
[X_1,[X_1,X_2]](x,y)=\left( \begin{array}{c}
0 \\ \gamma(\gamma-1) x^{\gamma-2}
\end{array} \right)\, ,\ \dots
$$
Thus, if $\gamma =1$, H\"ormander's condition is satisfied with $X_1$ and $[X_1,X_2]$.
In general, if $\gamma \ge 1$,  $\gamma $ iterated Lie brackets are required.

Theorem \ref{Main-thm} emphasizes that hypoellipticity is not sufficient for null controllability: Grushin's 
operator is hypoelliptic, but null controllability holds only when $\gamma =1$. 

The situation is similar for the Kolmogorov equation (\ref{Kolm}), where
$$X_0(x,v):=\left( \begin{array}{c}
v \\ 0
\end{array} \right)\, ,\quad
X_1(x,v):=\left( \begin{array}{c}
0 \\ 1
\end{array} \right)\, ,\quad
[X_1,X_2](x,v)=\left( \begin{array}{c}
1 \\ 0
\end{array} \right)\,.$$
Here again, null controllability holds 
and the first iterated Lie bracket is sufficient to satisfy H\"ormander's condition.

A general result which relates null controllability to the number of iterated Lie brackets that are necessary 
to satisfy H\"ormander's condition would be very interesting, but remains---for the time being---a
challenging open problem.

\subsection{Structure of the article}

Section \ref{Proof of Theorem OI} is devoted to the proof of Theorem \ref{thm:OI}.
In Subsection \ref{Well posedness} we recall useful results about the well-posedness of Grushin's equation. In Subsection~\ref{Fourier} we justify the Fourier decomposition 
of the solution to the adjoint system, which is needed for the proof of our main result.
In Subsection~\ref{Strategy} we present the strategy for the proof of Theorem \ref{thm:OI},
which relies on uniform observability estimates with respect to Fourier frequencies.
In Subsection \ref{Dissipation speed} we prove a preliminary result,
related to the dissipation rate of the Fourier components of the solution to the Grushin equation.
In Subsection \ref{subsec:Proof_positive}, we prove the positive statements of Theorem \ref{thm:OI},
thanks to an appropriate Carleman inequality.
In Subsection \ref{subsec:Proof_negative}, we show the negative statements of Theorem \ref{thm:OI},
thanks to appropriate test functions to falsify the observability inequality.
Then, in Section \ref{subsec:ccl}, we complete the proof of Theorem \ref{thm:OI}.
Finally, in Section \ref{sec:ccl}, we present several open problems and perspectives.

\section{Proof of Theorem \ref{thm:OI}}
\label{Proof of Theorem OI}

\subsection{Well posedness of the Cauchy-problem}
\label{Well posedness}

Let $H := L^2(\Omega;\mathbb{R})$, denote with $\langle\cdot,\cdot\rangle$ the scalar product in $H$ and by $\|\,\cdot\,\|_H := \langle\cdot,\cdot\rangle^{1/2}$ its norm. Define the scalar product
\begin{equation}
(f,g) := \int_\Omega\left(f_x g_x + |x|^{2\gamma} f_y g_y\right) dxdy
\end{equation}
for every $f$, $g$ in $C^\infty_0(\Omega)$, and set $V = \overline{C^\infty_0(\Omega)}^{|\,\cdot\,|_V}$, where $|f|_V := (f,f)^{1/2}$.

Observe that $H^1_0(\Omega)\subset V\subset H$, thus $V$ is dense in $H$.
We define the bilinear form $a$ on $V$ by
\begin{equation}
a(f,g) = - (f,g) \quad \forall f,\,g\in V\, .
\end{equation}
Moreover, set
\begin{equation}
D(A) = \left\{f\in V : \exists\, c>0 \text{ such that } |a(f,h)| \le c \|h\|_H \ \forall h\in V\right\}\, ,
\end{equation}
\begin{equation}\label{opA}
\langle Af,h\rangle = a(f,h)\qquad \forall h\in V\, .
\end{equation}

Then, we can apply a result by Lions \cite{lions68} (see also Theorem~1.18 in \cite{Zab}) and conclude that $(A,D(A))$ generates an analytic semigroup $S(t)$ of contractions on $H$. Note that $A$ is selfadjoint on $H$, and \eqref{opA} implies that
$$
Af = \partial_x^2 f + |x|^{2\gamma} \partial_y^2 f\qquad \text{almost everywhere in } \Omega\, .
$$

So, system (\ref{Grushin_CYpb}) can be recast in the form
\begin{equation}\label{eqabst}
\begin{cases}
f'(t) = Af(t) + u(t)\, \quad  t\in [0,T]\, , \\
f(0) = f_0\, ,
\end{cases}
\end{equation}
where $T>0$, $u\in L^2(0,T;H)$ and $f_{0} \in H$.

Let us now recall the definition of weak solutions to \eqref{eqabst}.
\begin{Def}[Weak solution]\label{Def:defweaksol}
Let $T>0$, $u\in L^2(0,T;H)$ and $f_{0} \in H$.
A function $f\in C([0,T];H)$ is a weak solution of \eqref{eqabst} if for every $h\in D(A)$ the function $\langle f(t),h\rangle$ is absolutely continuous on $[0,T]$ and for a.e. $t\in [0,T]$
\begin{equation}\label{eq:defweaksol}
\frac{d}{dt} \langle f(t),h\rangle = \langle f(t),A h\rangle + \langle u(t),h\rangle\, .
\end{equation}
\end{Def}
Note that, as showed in \cite{lions61}, condition \eqref{eq:defweaksol} is equivalent to the definition of solution by transposition, that is,
$$
\begin{array}{l}
\displaystyle \int\limits_{\Omega} 
[ f(t^{*},x,y)\varphi(t^{*},x,y) - f_{0}(x,y) \varphi(0,x,y) ]dx dy \\
\displaystyle  = \int\limits_{0}^{t^{*}} \int\limits_{\Omega} \left\{
f \left( 
\partial_t \varphi + \partial_x^2 \varphi + |x|^{2\gamma} \partial_y^2 \varphi \right)
+ u 1_{\omega} \varphi \right\} dx dy dt
\end{array}
$$
for every
$\varphi \in C^{2}([0,T] \times \Omega;\mathbb{R})$ and $t^{*} \in (0,T)$.

Let us recall that, for every $T>0$ and $u\in L^2(0,T;H)$, the mild solution $f\in C([0,T];H)$ of \eqref{eqabst} is defined as
\begin{equation}\label{varofconstfor}
f(t) = S(t)f_0 + \int_{0}^{t} S(t-s)u(s) ds \, ,\quad t\in [0,T]\, .
\end{equation}
From \cite{Ball}, we have that the mild solution to \eqref{eqabst} is also the unique weak solution in the sense of Definition \ref{Def:defweaksol}. The following existence and uniqueness result follows.

\begin{Prop} \label{Prop:existence}
For every $f_{0} \in H$, $T>0$ and 
$u \in L^{2}(0,T;H)$,
there exists a unique weak solution of the Cauchy problem (\ref{eqabst}).
This solution satisfies
\begin{equation}\label{C0continforsol}
\|f(t)\|_{H}
\leqslant   \|f_0\|_{H} + \|u\|_{L^2(0,T;H)} \quad\forall t \in [0,T]\,.
\end{equation}
Moreover, $f(t,\cdot) \in D(A)$ for a.e. $t\in (0,T)$.
\end{Prop}

\noindent \textbf{Proof:}
Inequality \eqref{C0continforsol} follows from \eqref{varofconstfor}. 
Moreover, since $S(\cdot)$ is analytic, $f(t) \in D(A)$ for a.e. $t\in (0,T)$.
 $\Box$

\subsection{Fourier decomposition}
\label{Fourier}

Let us consider the solution of (\ref{adjointL2}) in the sense of Definition \ref{Def:defweaksol}, that is, the solution of system (\ref{eqabst}) with $u = 0$.
The function $g$ belongs to $C([0,T];L^{2}(\Omega))$,
so $y \mapsto g(t,x,y)$ belongs to $L^{2}(0,1)$ for a.e. $(t,x) \in (0,T) \times (-1,1)$,
thus it can be developed in Fourier series in $y$
\begin{equation}\label{eq:g}
g(t,x,y)=\sum\limits_{n \in \mathbb{N}^*} g_n(t,x) \varphi_n(y)\, ,
\end{equation}
where
$$\varphi_n(y):=\sqrt{2} \sin(n\pi y) \quad\forall n \in \mathbb{N}^*$$
and
\begin{equation}\label{eq:gn}
g_n(t,x):=\int_0^1 g(t,x,y) \varphi_n(y) dy\quad \forall n \in \mathbb{N}^*\, .
\end{equation}

\begin{Prop}\label{Thm:wellposed}
For every $n\ge 1$, $g_n(t,x)$ is the unique weak solution of
\begin{equation} \label{adjointL2_n}
\left\lbrace \begin{array}{ll}
\partial_t g_n - \partial_x^2 g_n + (n\pi)^2 |x|^{2\gamma} g_n =0 & (t,x) \in (0,T)\times(-1,1)\,,\\
g_n(t,\pm 1)=0 &t\in(0,T)\,,\\
g_n(0,x)=g_{0,n}(x) &x\in (-1,1)\,.
\end{array}\right.
\end{equation}
\end{Prop}
For the proof we need the following characterization of the elements of $V$. We denote by $L^2_\gamma(\Omega)$ the space of all the square-integrable functions with respect to the measure $d\mu = |x|^{2\gamma}dxdy$.
\begin{Lem}\label{Lem:characV}
For every $g\in V$ there exist $\partial_x g\in L^2(\Omega)$, $\partial_y g\in L^2_\gamma(\Omega)$ such that
\begin{multline}
\int_\Omega \left(g(x,y) \partial_x\phi(x,y) + |x|^{2\gamma} g(x,y) \partial_y\phi(x,y)\right) dxdy \\
 = -\int_\Omega \left(\partial_x g(x,y) + |x|^{2\gamma}\partial_y g(x,y)\right)\phi(x,y) dxdy
\end{multline}
for every $\phi\in C^\infty_0(\Omega)$.
\end{Lem}

\noindent \textbf{Proof:} Let $g\in V$, and consider a sequence $(g^n)_{n\ge 1}$ in $C^\infty_0(\Omega)$ such that $g^n\to g$ in $V$, that is
$$
\int_\Omega\left[(g^n - g)_x^2 + |x|^{2\gamma}(g^n - g)_y^2\right]dxdy\to 0\quad \text{as $n\to +\infty$}\, .
$$
Thus, $(\partial_x g^n)_{n\ge 1}$ is a Cauchy sequence in $L^2(\Omega)$, and $(\partial_y g^n)_{n\ge 1}$ is a Cauchy sequence in $L^2(\Omega,|x|^{2\gamma}dxdy)$, so there exist $h\in L^2(\Omega)$ and $k\in L^2(\Omega, |x|^{2\gamma}dxdy)$ such that $\partial_x g^n\to h$ in $L^2(\Omega)$ and $\partial_y g^n\to k$ in $L^2(\Omega, |x|^{2\gamma}dxdy)$. Hence,
\begin{displaymath}
\begin{array}{ccc}
\displaystyle \int_\Omega \left(g^n \partial_x \phi + |x|^{2\gamma}g^n\partial_y \phi\right)dxdy & = & \displaystyle -\int_\Omega\left(\partial_x g^n \phi + |x|^{2\gamma}\partial_y g^n\phi\right)dxdy \\
\displaystyle \Big\downarrow &  & \displaystyle \Big\downarrow \\
\displaystyle \int_\Omega\left(g \partial_x\phi + |x|^{2\gamma}g\partial_y\phi\right)dxdy & = & \displaystyle -\int_\Omega\left(h \phi + |x|^{2\gamma}k\phi\right)dxdy
\end{array}
\end{displaymath}
as $n\to +\infty$. This yields the conclusion with $\partial_x g = h$ and $\partial_y g = k$. $\Box$
\smallskip
\smallskip

For any $n\ge 1$, system \eqref{adjointL2_n} is a first order Cauchy problem, that admits a unique weak solution
$$
\tilde{g}_n\in C^1((0,T);L^2(-1,1))\cap C([0,T];H^1_0(-1,1))\cap L^2(0,T;H^2(-1,1))
$$
which satisfies
\begin{multline}\label{eq:weakdefgn}
\frac{d}{dt}\left(\int_{-1}^1 \tilde{g}_n(t,x) \psi(x) dx\right) \\
+ \int_{-1}^1 \Big[\tilde{g}_{n,x}(t,x)\psi_x(x) + (n\pi)^2 |x|^{2\gamma} \tilde{g}_n(t,x) \psi(x)\Big]dx = 0
\end{multline}
for every $\psi\in H^1_0(-1,1)$.
\smallskip
\smallskip

\noindent \textbf{Proof of Proposition \ref{Thm:wellposed}:} In order to verify that the $n$th Fourier coefficient of $g$, defined by \eqref{eq:gn}, satisfies system \eqref{adjointL2_n}, observe that
$$
g_n(0,.) = g_{0,n}(.)\, ,\qquad g_n(t,\pm 1) = 0\quad \forall t\in (0,T)
$$
and
$$
g_n(t,x) \in C^1((0,T);L^2(-1,1))\cap C([0,T];H^1_0(-1,1))\, .
$$
Thus, it is sufficient to prove that $g_n$ fulfills condition \eqref{eq:weakdefgn}. Indeed, using the identity (\ref{eq:gn}), for all $\psi\in H^1_0(-1,1)$,
\begin{multline}\label{eq:step1}
\frac{d}{dt}\left(\int_{-1}^1 g_n \psi dx\right) + \int_{-1}^1 \left(g_{n,x}\psi_x + (n\pi)^2 |x|^{2\gamma}g_n \psi\right)dx \\
= \int_{-1}^1 \int_0^1\left\{g_t\varphi_n\psi + g_x\varphi_n\psi_x + (n\pi)^2|x|^{2\gamma} g\varphi_n\psi\right\}dydx\, .
\end{multline}
On the other hand, choosing $h(x,y) = \psi(x)\varphi_n(y)\in V$ in \eqref{eq:defweaksol},
\begin{multline}\label{eq:step2}
0 = \int_0^1\int_{-1}^1(g_t - Ag)\psi\varphi_ndxdy \\
= \int_0^1\int_{-1}^1g_t\psi\varphi_n dxdy + \int_0^1\int_{-1}^1 \left(g_x\psi_x\varphi_n + |x|^{2\gamma}g_y \psi\varphi_{n,y}\right)dxdy \\
= \int_0^1\int_{-1}^1g_t\psi\varphi_n dxdy + \int_0^1\int_{-1}^1\left( g_x\psi_x\varphi_n + (n\pi)^2 |x|^{2\gamma}g \psi\varphi_{n}\right)dxdy\, ,
\end{multline}
where (in the last identity) we have used Lemma \ref{Lem:characV}. Combining \eqref{eq:step1} and \eqref{eq:step2} completes the proof. $\Box$

\subsection{Strategy for the proof of Theorem \ref{thm:OI}}
\label{Strategy}

\noindent 
Let $g$ be the solution of (\ref{adjointL2}). Then, $g$ can be represented as in \eqref{eq:g}, and we emphasize that, for a.e. $t\in (0,T)$, and for every $-1\leqslant a_1 <b_1 \leqslant 1$,
$$\int_{(a_1,b_1) \times (0,1)} |g(t,x,y)|^2 dx dy = \sum\limits_{n=1}^{\infty} \int_{a_1}^{b_1} |g_n(t,x)|^2 dx$$
(Bessel-Parseval equality). Thus, in order to prove Theorem \ref{thm:OI}, it is sufficient to study the 
observability of system (\ref{adjointL2_n}) uniformly with respect to $n \in \mathbb{N}^*$.

\begin{Def}[Uniform observability] Let $0<a<b<1$ and $T>0$. System (\ref{adjointL2_n}) is {\em observable in $(a,b)$ in time $T$ uniformly with respect to $n \in \mathbb{N}^*$} if
there exists $C>0$ such that,
for every $n \in \mathbb{N}^*$, $g_{0,n} \in L^2(-1,1)$, 
the solution of (\ref{adjointL2_n}) satisfies
$$\int_{-1}^1 |g_n(T,x)|^2 dx \leqslant C \int_0^T \int_a^b |g_n(t,x)|^2 dx\, .$$

System (\ref{adjointL2_n}) is {\em observable in $(a,b)$  uniformly with respect to $n \in \mathbb{N}^*$} if
there exists $T>0$ such that it is observable in $(a,b)$ in time $T$ uniformly with respect to $n \in \mathbb{N}^*$.
\end{Def}

Theorem \ref{thm:OI} is a consequence of the following statement.

\begin{thm} \label{thm:OI_n}
We assume $0<a<b<1$.
\begin{enumerate}
\item If $\gamma \in (0,1)$, then system (\ref{adjointL2_n}) is observable in $(a,b)$ in any time $T>0$ uniformly with respect to $n \in \mathbb{N}^*$.
\item If $\gamma=1$, there exists $T^*>0$ such that
\begin{itemize}
\item for every $T>T^*$, system (\ref{adjointL2_n}) is observable in $(a,b)$  in time $T$ uniformly with respect to $n \in \mathbb{N}^*$,
\item for every $T<T^*$, system (\ref{adjointL2_n}) is not observable in $(a,b)$  in time $T$ uniformly with respect to $n \in \mathbb{N}^*$.
\end{itemize}
\item If $\gamma>1$, system (\ref{adjointL2_n}) is not observable in $(a,b)$  uniformly with respect to $n \in \mathbb{N}^*$.
\end{enumerate}
\end{thm}

The strategy of the proof for the positive statements of Theorem \ref{thm:OI_n} is standard and relies on two key ingredients:
\begin{itemize}
\item an explicit decay rate for the solutions of (\ref{adjointL2_n}),
\item a favorable estimate for the observability constant associated to the equation (\ref{adjointL2_n}) 
and the observation domain $(a,b)$.
\end{itemize}
This strategy has already been used in \cite{JMC-Guerrero}, \cite{Kolm}.
The proof of the negative statements of Theorem \ref{thm:OI_n} relies on the use of appropriate test functions
that falsify uniform observability.

Let us recall that explicit bounds on the observability constant of  the heat equation with a potential are already known.
\begin{thm}
Let $-1<a<b<1$. There exists $c>0$ such that,
for every $T>0$, $\alpha$, $\beta \in L^{\infty}((0,T) \times (-1,1))$, $g_0 \in L^2(-1,1)$,
the solution of 
$$\left\lbrace \begin{array}{ll}
\partial_t g - \partial_x^2 g + \beta \partial_x g  + \alpha g = 0 &  (t,x) \in [0,T]\times (-1,1)\, ,
\\
g(t,\pm 1)=0 & t\in [0,T]\, ,
\\
g(0,x)=g_{0}(x) & x\in (-1,1)\, ,
\end{array}\right.$$
satisfies
$$\int_{-1}^{1} |g(T,x)|^2 dx \leqslant e^{c H(T,\|\alpha\|_{\infty} , \|\beta\|_{\infty})}
\int_0^T \int_a^b |g(t,x)|^2 dx dt\, ,$$
where $H(T,A,B):=1 + \frac{1}{T} + T A + A^{2/3} + (1+T)B^2$.
\end{thm}
For the proof of this result, we refer to 
\cite[Theorem 1.3]{Cara-Zuazua-ChaleurLineaire} in the case $\beta \equiv 0$ 
and to \cite[Theorem 2.3]{Doubova-EFC-MGB-EZ} in the case $\beta \neq 0$.
The optimality of the power $2/3$ of $A$ in $H(T,A,B)$ has been proved in \cite{Duychaert-Zhang-Zuazua}. 

The positive statement of Theorem \ref{thm:OI_n} may be seen 
as an improvement of the above estimate (relatively to the asymptotic behavior when $n \rightarrow + \infty$), 
in the particular case of equation (\ref{adjointL2_n}).

\subsection{Dissipation speed}
\label{Dissipation speed}

Let us introduce, for every $n \in \mathbb{N}^*, \gamma>0$, the operator $A_{n,\gamma}$ defined by
\begin{equation} \label{def:An}
\begin{array}{ll}
D(A_{n,\gamma}):=H^2 \cap H^1_0(-1,1;\mathbb{R})\,, 
&
A_{n,\gamma} \varphi := -\varphi'' + (n\pi)^2 |x|^{2\gamma} \varphi\, .
\end{array}
\end{equation}
The smallest eigenvalue of $A_{n,\gamma}$ is given by
\begin{equation}\label{charfirsteigenvalue}
\displaystyle\lambda_{n,\gamma} = \min \left\{
\frac{ \int_{-1}^1 \left[ v'(x)^2 + (n\pi)^2 |x|^{2\gamma} v(x)^2 \right] dx }{\int_{-1}^1 v(x)^2 dx } ;
v \in H^1_0(-1,1) \right\}\, .
\end{equation}
We are interested in the asymptotic behavior (as $n \rightarrow + \infty$) of $\lambda_{n,\gamma}$,
which quantifies the dissipation speed of the solution of (\ref{adjointL2_n}).

Thanks to a simple heuristic computation, one may expect that, for every $\gamma>0$,  
$\lambda_{n,\gamma}$ behaves like  $C(\gamma) n^{\frac{2}{1+\gamma}}$.
Indeed, if we consider the eigenvector $v_{n,\gamma}$ 
\begin{equation}\label{systeigenfunc}
\left\lbrace \begin{array}{l}
- v_{n,\gamma}''(x)+(n\pi)^2 |x|^{2\gamma} v_{n,\gamma}(x)=\lambda_{n,\gamma} v_{n,\gamma}(x) \quad x\in (-1,1)\,,\\
v_{n,\gamma}(\pm 1)=0\, ,
\end{array}\right.
\end{equation}
and the change of variable $y:=l x$, $l=(n\pi)^{1/(1+\gamma)}$, $v_{n,\gamma}(x)=\psi_{n,\gamma}(y)$, we get
$$\left\lbrace \begin{array}{l}
- \psi_{n,\gamma}''(y)+|y|^{2\gamma} \psi_{n,\gamma}(y)=\lambda_{n,\gamma}  (n\pi)^{\frac{-2}{1+\gamma}} \psi_{n,\gamma}(y)\, \quad y\in (-l, l)\, , \\
\psi_{n,\gamma}(\pm l )=0\, .
\end{array}\right.$$
In order to prove two results related to this conjecture, we need the following lemma.
\begin{Lem}\label{Lem:firsteigeven}
Problem \eqref{systeigenfunc} admits a unique positive solution with norm one. Moreover, $v_{n,\gamma}$ is even.
\end{Lem}

\noindent \textbf{Proof:} Since \eqref{systeigenfunc} is a Sturm-Liouville problem, it is well-known that its first eigenvalue is simple, and the associated eigenfunction has no zeros. Thus, we can choose $v_{n,\gamma}$ to be strictly positive everywhere. Moreover, by normalization, we can find a unique positive solution satisfying the condition $\|v_{n,\gamma}\|_{L^2(-1,1)} = 1$. Finally, $v_{n,\gamma}$ is even. Indeed, if not so, let us consider the function $w(x) = v_{n,\gamma}(|x|)$. Then, $w$ still belongs to $H^1_0(-1,1)$, it is a weak solution of \eqref{systeigenfunc} and it does not increase the functional in (\ref{charfirsteigenvalue}), i.e.
$$
\frac{ \int_{-1}^1 \left[ w'(x)^2 + (n\pi)^2 |x|^{2\gamma} w(x)^2 \right] dx }{\int_{-1}^1 w(x)^2 dx }\le 
\frac{ \int_{-1}^1 \left[ v_{n,\gamma}'(x)^2 + (n\pi)^2 |x|^{2\gamma} v_{n,\gamma}(x)^2 \right] dx }{\int_{-1}^1 v_{n,\gamma}(x)^2 dx }\, .
$$
The coefficients of the equation in \eqref{systeigenfunc} being regular, we deduce that $w$ is a classical solution of \eqref{systeigenfunc}. Since $\lambda_{n,\gamma}$ is simple, it follows $v_{n,\gamma}(x) = v_{n,\gamma}(|x|)$.
$\Box$

The following result turns out to be the key point of the proof of Theorem~\ref{thm:OI_n}.

\begin{Prop} \label{Prop:1st_eigenvalue}
\begin{enumerate}
\item For every $\gamma \in (0,1]$, there exists $c_*=c_*(\gamma)>0$ such that
$$\lambda_{n,\gamma} \geqslant c_* n^{\frac{2}{1+\gamma}}\qquad \forall\, n \in \mathbb{N}^*\,.$$
\item For every $\gamma>0$, there exists $c^*=c^*(\gamma)>0$ such that
$$\lambda_{n,\gamma} \leqslant c^* n^{\frac{2}{1+\gamma}}\qquad \forall\, n \in \mathbb{N}^*\,.$$
\end{enumerate}
\end{Prop}

\noindent \textbf{Proof:} 
First, taking $\gamma \in (0,1]$, let us prove the first part of the conclusion.
Thanks to Lemma \ref{Lem:firsteigeven}, we have
$$\lambda_{n,\gamma} = \min \left\{
\frac{ \int_0^1 \left[ v'(x)^2 + (n\pi)^2 x^{2\gamma} v(x)^2 \right]dx }{ \int_0^1 v(x)^2 dx } ;
v \in H^1_0(-1,1) \right\}.$$
Thus, our goal is to prove the existence of $c_0=c_0(\gamma)>0$ such that
$$\int_0^1 \left[ v'(x)^2 + (n\pi)^2 x^{2\gamma} v(x)^2 \right] dx \geqslant c_0  (n\pi)^{\frac{2}{1+\gamma}} \int_0^1 v(x)^2 dx \quad
\forall v \in H^1_0(-1,1)\,,$$
or, equivalently, the existence of $c_1=c_1(\gamma)>0$ such that
\begin{equation} \label{goal_1}
\int_0^1 v(x)^2 dx \leqslant c_1 \int_0^1 \left[ \alpha_n^2 v'(x)^2 + \left( \frac{x}{\alpha_n} \right)^{2\gamma} v(x)^2 \right] dx \quad
\forall v \in H^1_0(-1,1)\,,
\end{equation}
where $\alpha_n:=(n\pi)^{\frac{-1}{1+\gamma}}$.
First, let us emphasize that
$$\int_{\alpha_n}^1 v(x)^2 dx \leqslant \int_{\alpha_n}^{1} \left( \frac{x}{\alpha_n} \right)^{2\gamma} v(x)^2 dx\, .$$
Thus, in order to prove (\ref{goal_1}), it is sufficient to find $c_2=c_2(\gamma)>0$ such that
\begin{equation} \label{goal_2}
c_2 \int_0^{\alpha_n} v(x)^2 dx \leqslant  \int_0^1 \left[ \alpha_n^2 v'(x)^2 + \left( \frac{x}{\alpha_n} \right)^{2\gamma} v(x)^2 \right] dx\quad
\forall v \in H^1_0(-1,1)\,.
\end{equation}

\emph{First step: Let us prove that, for all $v \in H^1_0(-1,1)$,}
\begin{equation} \label{goal_3}
\gamma \int_0^{\alpha_n} v(x)^2 dx \leqslant  
\alpha_n v(\alpha_n)^2 + \int_0^1 \left[ \alpha_n^2 v'(x)^2 + \left( \frac{x}{\alpha_n} \right)^{2\gamma} v(x)^2 \right] dx\, .
\end{equation}
Let $v \in H^1_0(-1,1)$. We have
\begin{multline*}
\int_0^1 \left[ \alpha_n^2 v'(x)^2 + \left( \frac{x}{\alpha_n} \right)^{2\gamma} v(x)^2 \right]dx \geqslant
\int_0^{\alpha_n} \left[ \alpha_n^2 v'(x)^2 + \left( \frac{x}{\alpha_n} \right)^{2\gamma} v(x)^2 \right] dx
\\ \geqslant
-2 \int_0^{\alpha_n} \alpha_n v'(x) \left( \frac{x}{\alpha_n} \right)^{\gamma} v(x) dx =
- \alpha_n v(\alpha_n)^2 + \gamma \int_0^{\alpha_n} \left( \frac{\alpha_n}{x} \right)^{1-\gamma} v(x)^2 dx
\\ \geqslant
- \alpha_n v(\alpha_n)^2 + \gamma \int_0^{\alpha_n}  v(x)^2 dx
\end{multline*}
because $1-\gamma \geqslant 0$. This proves inequality (\ref{goal_3}).
\smallskip
\smallskip

\emph{Second step: Let us prove the existence of $c_3=c_3(\gamma)>0$ such that}
\begin{equation} \label{goal_4}
 \alpha_n v(\alpha_n)^2
\leqslant
c_3 \int_0^1 \left[ \alpha_n^2 v'(x)^2 + \left( \frac{x}{\alpha_n} \right)^{2\gamma} v(x)^2 \right] dx\quad
\forall v \in H^1_0(-1,1)\,.
\end{equation}
Let $v \in H^1_0(-1,1)$. Since
$$v(\alpha_n)=v(x)+\int_x^{\alpha_n} v'(s) ds\quad \forall x \in (0,\alpha_n)\,,$$
we have
$$\int_0^{\alpha_n} x^{2\gamma} v(\alpha_n)^2 dx
\leqslant
2 \int_0^{\alpha_n} x^{2\gamma} v(x)^2 dx 
+
2 \int_0^{\alpha_n} x^{2\gamma} (\alpha_n-x) dx \int_{0}^{\alpha_n} v'(s)^2 ds\,,$$
which implies
$$\frac{\alpha_n^{2\gamma+1}}{2\gamma+1} v(\alpha_n)^2
\leqslant
2 \int_0^{\alpha_n} x^{2\gamma} v(x)^2 dx 
+2 \left( \alpha_n \frac{\alpha_n^{2\gamma+1}}{2\gamma+1} - \frac{\alpha_n^{2\gamma+2}}{2\gamma+2} \right) \int_0^{\alpha_n} v'(x)^2 dx\,.$$
Multiplying both sides by $(2\gamma+1)\alpha_{n}^{-2\gamma}$, we deduce
$$\alpha_n v(\alpha_n)^2 
\leqslant
2(2\gamma+1)\int_0^{\alpha_n} \left( \frac{x}{\alpha_n} \right)^{2\gamma} v(x)^2 dx 
+
\frac{1}{\gamma+1} \int_0^{\alpha_n} \alpha_n^2 v'(x)^2 dx\,.$$
Hence, (\ref{goal_4}) holds with $c_3(\gamma)=2(2\gamma+1)$. Combining (\ref{goal_3}) and (\ref{goal_4}) gives (\ref{goal_2}) with $c_2=\gamma/(1+c_3)$.

Now, let $\gamma >0$ and let us prove the second statement of Proposition \ref{Prop:1st_eigenvalue}.
For every $k>1$ we consider the function
$\varphi_k(x):= (1 - k|x|)^+$, that belongs to $H^1_0(-1,1)$. Easy computations show that
$$\begin{array}{c}
\displaystyle \int_{-1}^1 \varphi_k(x)^2 dx = \frac{2}{3k}\, ,\ \displaystyle\int_{-1}^1 \varphi_k'(x)^2 dx = 2k\, , \
\displaystyle \int_{-1}^1 |x|^{2\gamma} \varphi_k(x)^2 dx = 2c(\gamma) k^{-1-2\gamma}\, ,
\end{array}$$
where
$$c(\gamma) := \left( \frac{1}{2\gamma+1} - \frac{1}{\gamma+1} + \frac{1}{2\gamma+3} \right).$$
Thus, $\lambda_{n,\gamma} \leqslant f_{n,\gamma}(k):=3[k^2 + (\pi n)^2 c(\gamma) k^{-2\gamma}]$ for all $k>1$.
Since $f_{n,\gamma}$ attains its minimum at $\bar{k} = \tilde{c}(\gamma) n^{\frac{1}{\gamma+1}}$, we have $\lambda_{n,\gamma} \leqslant f_{n,\gamma}(\bar{k})=C(\gamma) n^{\frac{2}{\gamma+1}}$.$\Box$

\subsection{Proof of the positive statements of Theorem \ref{thm:OI_n}}
\label{subsec:Proof_positive}

The goal of this section is the proof of the following results:
\begin{itemize}
\item if $\gamma \in (0,1)$ and $T>0$, then system (\ref{adjointL2_n}) is
observable in time $T$ uniformly with respect to $n$;
\item if $\gamma=1$, there exists $T_1>0$ such that, for every $T>T_1$, system (\ref{adjointL2_n}) is observable in time $T$ uniformly with respect to $n$.
\end{itemize}
The proof of these results relies on a new Carleman estimate for the solutions of (\ref{adjointL2_n}). 

Let $\gamma \in (0,1]$, $T>0$, and fix $n \in \mathbb{N}^*$ all over the proof. In order to simplify the notation, we write $g$ and $g_0$ instead of $g_n$ and $g_{0,n}$. Let $a', b'$ be such that $a < a' < b' < b$.
\smallskip
\smallskip

\textbf{\underline{First case: $\gamma \in [1/2,1]$}}
\smallskip

In order to deduce the Carleman inequality, we define a weight function
\begin{equation} \label{def-alpha}
\alpha(t,x):=\frac{M \beta (x)}{t(T-t)}\, ,\quad (t,x) \in (0,T) \times \mathbb{R}\,,
\end{equation}
where $\beta \in C^{2}(\mathbb{R};\mathbb{R}_{+})$ satisfies
\begin{equation} \label{hyp-beta}
\beta \geqslant 1 \text{  on  } (-1,1)\,,
\end{equation}
\begin{equation} \label{hyp-beta'-(a,b)}
|\beta'|>0 \text{  on  } [-1,a'] \cup [b',1]\,,
\end{equation}
\begin{equation} \label{hyp_beta'-bord}
\beta'(1)>0\,,\quad \beta'(-1)<0\,,
\end{equation}
\begin{equation} \label{hyp-beta''}
\beta'' <0 \text{ on  }  [-1,a'] \cup [b',1]\,,
\end{equation}
and $M>0$ will be chosen later on. We also introduce the function
\begin{equation} \label{def-z}
z(t,x):=g(t,x) e^{-\alpha(t,x)}\, ,
\end{equation}
that satisfies
\begin{equation} \label{P123}
P_{1}z + P_{2} z = P_{3} z\, ,
\end{equation}
where
\begin{equation} \label{def:P1-P2-P3}
\begin{array}{c}
P_{1}z:= - \frac{\partial^{2} z}{\partial x^{2}} +(\alpha_{t}-\alpha_{x}^{2}) z + (n\pi)^2 |x|^{2\gamma} z \, ,\quad
P_{2}z:= \frac{\partial z}{\partial t} - 2 \alpha_{x} \frac{\partial z}{\partial x}\, , \\
P_{3}z:= \alpha_{xx} z\,.
\end{array}
\end{equation}
We develop the classical proof, taking the $L^{2}(Q)$-norm in the identity (\ref{P123}),
then developing the double product, which leads to
\begin{equation} \label{P1P2<P3}
\int_{Q} P_{1}z P_{2}z  \leqslant \frac{1}{2} \int_{Q} |P_{3}z|^{2}\, ,
\end{equation}
where $Q:=(0,T) \times (-1,1)$ and we compute precisely each term.
\smallskip

\textbf{Terms concerning $- \partial_x^{2} z$: }
Integrating by parts, we get
\begin{equation} \label{carl1}
- \int_{Q} \frac{\partial^{2} z}{ \partial x^{2}} \frac{\partial z}{\partial t} dxdt
=  \int_{Q} 
\frac{\partial z}{ \partial x} \frac{\partial^{2} z}{ \partial t \partial x}dxdt
= \int_{0}^{T} \frac{1}{2} \frac{d}{dt} \int_{-1}^1 
\Big| \frac{\partial z}{\partial x} \Big|^{2} dxdt
=0\,,
\end{equation}
because $\partial_t z(t,\pm 1)=0$ and $z(0) \equiv z(T) \equiv 0$, which is a consequence of assumptions
(\ref{def-z}), (\ref{def-alpha}) and (\ref{hyp-beta}).
Moreover,
\begin{multline} \label{carl2}
\int_{Q}  \frac{\partial^{2} z}{ \partial x^{2}} 
2 \alpha_{x} \frac{\partial z}{\partial x} dxdt
= 
- \int_{Q} \Big| \frac{\partial z}{ \partial x} \Big|^{2} \alpha_{xx} dxdt
\\
+ \int_{0}^{T} \left(
\alpha_{x}(t,1) \Big|\frac{\partial z}{ \partial x}(t,1) \Big|^{2} -
\alpha_{x}(t,-1)    \Big|\frac{\partial z}{ \partial x}(t,-1)    \Big|^{2} \right) dt 
\, .
\end{multline}

\textbf{Terms concerning $(\alpha_{t}-\alpha_{x}^{2}) z$: }
Again integrating by parts, we have
\begin{equation} \label{carl3}
\int_{Q} (\alpha_{t}-\alpha_{x}^{2}) z \frac{\partial z}{\partial t} dxdt =
-\frac{1}{2} \int_{Q} (\alpha_{t}-\alpha_{x}^{2})_{t} |z|^{2} dxdt\, .
\end{equation}
Indeed, the boundary terms at $t=0$ and $t=T$ vanish because, thanks to 
(\ref{def-z}), (\ref{def-alpha}), (\ref{hyp-beta}),
$$|(\alpha_{t}-\alpha_{x}^{2}) |z|^{2}|
 \leqslant 
\frac{1}{[t(T-t)]^{2}} e^{\frac{-M}{t(T-t)}}
| M(T-2t)\beta + (M \beta')^{2} |\cdot |g|^{2}$$
tends to zero when $t \rightarrow 0$ and $t \rightarrow T$,
for every $x \in [-1,1]$.
Moreover,
\begin{equation} \label{carl4}
- 2 \int_{Q}  (\alpha_{t}-\alpha_{x}^{2}) z \alpha_{x} \frac{\partial z}{\partial x} dxdt =
\int_{Q} [ (\alpha_{t}-\alpha_{x}^{2}) \alpha_{x} ]_{x} |z|^{2} dxdt\, ,
\end{equation}
thanks to an integration by parts in the space variable.
\smallskip

\textbf{Terms concerning $(n\pi)^2 |x|^{2\gamma} z$: }
First, since $z(0) \equiv z(T) \equiv 0$,
\begin{equation} \label{carl5}
\int_Q (n\pi)^2 |x|^{2\gamma} z \frac{\partial z}{\partial t} dxdt
=
\frac{1}{2} \int_0^T \frac{d}{dt} \int_{-1}^{1} (n\pi)^2 |x|^{2\gamma} |z|^2 dx dt
= 0\, .
\end{equation}
Furthermore, thanks to an integration by parts in the space variable,
\begin{equation} \label{carl6}
-2 \int_Q (n\pi)^2 |x|^{2\gamma} z \alpha_x \frac{\partial z}{\partial x} dxdt
=
\int_Q  [ n^2 \pi^2 |x|^{2\gamma} \alpha_x ]_x z^2 dxdt\, .
\end{equation}
Combining (\ref{P1P2<P3}), (\ref{carl1}), (\ref{carl2}), (\ref{carl3}),
(\ref{carl4}), (\ref{carl5}) and (\ref{carl6}), we conclude that
\begin{multline} \label{In-Carl-1}
\displaystyle\int_{Q} |z|^{2} \Big\{ 
- \frac{1}{2} (\alpha_{t}-\alpha_{x}^{2})_{t} 
+  [ (\alpha_{t}-\alpha_{x}^{2}) \alpha_{x} ]_{x}
+n^2 \pi^2 [ |x|^{2\gamma} \alpha_x ]_x 
- \frac{1}{2}\alpha_{xx}^{2}
\Big\} dxdt
\\
+ \int_{0}^{T} \left(
\alpha_{x}(t,1) \Big|\frac{\partial z}{ \partial x}(t,1) \Big|^{2} -
\alpha_{x}(t,-1)    \Big|\frac{\partial z}{ \partial x}(t,-1)    \Big|^{2} \right) dt 
\\
\displaystyle - \int_{Q} \Big| \frac{\partial z}{\partial x} \Big|^{2} \alpha_{xx} dxdt
\leqslant 0\,.
\end{multline}
In view of (\ref{hyp_beta'-bord}), we have $\alpha_{x}(t,1) \geqslant 0$ and
$\alpha_{x}(t,-1) \leqslant 0$, thus (\ref{In-Carl-1}) yields
\begin{multline} \label{In-Carl-1-BIS}
\int_{Q} |z|^{2} \Big\{ 
- \frac{1}{2} (\alpha_{t}-\alpha_{x}^{2})_{t} 
+  [ (\alpha_{t}-\alpha_{x}^{2}) \alpha_{x} ]_{x}
- \frac{1}{2}\alpha_{xx}^{2}
+n^2 \pi^2 [ |x|^{2\gamma} \alpha_x ]_x 
\Big\} dxdt
\\
- \int_{Q} \Big| \frac{\partial z}{\partial x} \Big|^{2} \alpha_{xx} dxdt
\leqslant 0\, .
\end{multline}
Now, in the left hand side of (\ref{In-Carl-1-BIS}) we separate the terms on
$(0,T) \times (a',b')$ and those on
$(0,T)\times [(-1,a') \cup (b',1)]$. One has
\begin{equation} \label{borne-dz}
\begin{array}{c}
\displaystyle -\alpha_{xx}(t,x) \geqslant \frac{C_{1} M}{t(T-t)} \qquad
\forall x \in [-1,a'] \cup [b',1]\,, \ t \in (0,T)\,,
\\
\displaystyle |\alpha_{xx}(t,x)| \leqslant \frac{C_{2} M}{t(T-t)} \qquad
\forall x \in [a',b']\,, \ t \in (0,T)\,,
\end{array}
\end{equation}
where $C_{1}:=\min \{ -\beta''(x) ; x \in [-1,a'] \cup [b',1] \}$ 
is positive thanks to the assumption (\ref{hyp-beta''}) and
$C_{2}:=\sup\{ |\beta''(x)| ; x \in [a',b'] \}$.
Moreover,
$$\begin{array}{l}
\displaystyle -\frac{1}{2}(\alpha_{t}-\alpha_{x}^{2})_{t} 
+  [ (\alpha_{t}-\alpha_{x}^{2}) \alpha_{x} ]_{x} 
- \frac{1}{2} \alpha_{xx}^{2} = \frac{1}{(t(T-t))^3} \Big\{M\beta (3Tt - T^{2} -3t^{2}) \\
\qquad \displaystyle + M^{2} \Big[ (2t-T) \beta'' \beta - \frac{t(T-t)\beta''^2}{2}  \Big]
- 3M^{3} \beta''  \beta '^{2}
 \Big\}\, .
\end{array}$$
Hence, owing to (\ref{hyp-beta'-(a,b)}) and (\ref{hyp-beta''}),
there exist $M_{1}=M_{1}(T,\beta)>0$, $C_{3}=C_{3}(\beta)>0$ and $C_{4}=C_{4}(T,\beta)>0$ such that, for every $M \geqslant M_{1}$ and $t \in (0,T)$,
\begin{multline} \label{borne-z}
\displaystyle -\frac{1}{2}(\alpha_{t}-\alpha_{x}^{2})_{t} 
\displaystyle +  [ (\alpha_{t}-\alpha_{x}^{2}) \alpha_{x} ]_{x} 
\displaystyle - \frac{1}{2} \alpha_{xx}^{2}
\geqslant \frac{C_{3}M^{3}}{[t(T-t)]^{3}}\quad
\forall x \in [-1,a'] \cup [b',1]\,, 
\\ 
\displaystyle \Big| -\frac{1}{2}(\alpha_{t}-\alpha_{x}^{2})_{t} 
+  [ (\alpha_{t}-\alpha_{x}^{2}) \alpha_{x} ]_{x} 
- \frac{1}{2} \alpha_{xx}^{2} \Big|
\leqslant \frac{C_{4}M^{3}}{[t(T-t)]^{3}}\quad
\forall x \in [a',b']\,.
\end{multline}

Using (\ref{In-Carl-1-BIS}), (\ref{borne-dz}) and (\ref{borne-z}), we deduce, for every $M \geqslant M_{1}$,
\begin{multline} \label{In-Carl-2}
\displaystyle\int_{0}^{T} \int\limits_{(-1,a')\cup(b',1)}
\frac{C_{1} M}{t(T-t)} \Big| \frac{\partial z}{\partial x} \Big|^{2} dxdt
\\
+ \displaystyle\int_{0}^{T} \int\limits_{(-1,a')\cup(b',1)}
\left[
\frac{C_{3} M^{3}}{(t(T-t))^{3}} |z|^{2}
+ (n\pi)^2 [ |x|^{2\gamma} \alpha_x ]_x |z|^2
\right]
dxdt
\\
\displaystyle\leqslant 
\int_{0}^{T} \int_{a'}^{b'}
\left[
\frac{C_{2} M}{t(T-t)}  \Big| \frac{\partial z}{\partial x} \Big|^{2}
+ \frac{C_{4} M^{3}}{(t(T-t))^{3}} |z|^{2}
- (n\pi)^2 [ |x|^{2\gamma} \alpha_x ]_x |z|^2
\right] dxdt\, .
\end{multline}
Moreover, for every $x \in (-1,1)$, we have
\begin{multline*}
\displaystyle|  (n\pi)^2 [ |x|^{2\gamma} \alpha_x ]_x |
=
\displaystyle \frac{M(n\pi)^2}{t(T-t)} \Big| 2\gamma \text{sign}(x) |x|^{2\gamma-1} \beta'(x) + |x|^{2\gamma} \beta''(x) \Big|
\leqslant
\frac{C_5 n^2 M}{t(T-t)}\, ,
\end{multline*}
where
$C_5:= \pi^2 \max\{ 2 \gamma |x|^{2\gamma-1} |\beta'(x)| + |x|^{2\gamma} |\beta''(x) |;x \in [-1,1] \}$
is finite because $2\gamma-1 \geqslant 0$. 
From now on, we take 
\begin{equation} \label{def:M}
M = M(T,\beta,n):=\max \{ 1 , M_{1}(T,\beta) , M_{2}(T,\beta,n) \}\, ,
\end{equation}
where $M_2=M_2(T,n)$ is defined by
\begin{equation} \label{def:M2}
M_2:=\sqrt{\frac{2C_5}{C_3}} n \left( \frac{T}{2} \right)^{2}\,.
\end{equation}
Since
$$|  (n\pi)^2 [ |x|^{2\gamma} \alpha_x ]_x |
\leqslant \frac{C_3 M^3}{2[t(T-t)]^3}\quad \forall (t,x) \in Q\,,$$
we conclude that
\begin{multline} \label{In-Carl-3}
\int_{0}^{T} \int\limits_{(-1,a') \cup (b',1)}  
\frac{C_{3} M^{3}}{2 (t(T-t))^{3}} |z|^{2} dxdt
\\
\leqslant 
\int_{0}^{T} \int_{a'}^{b'} 
\frac{C_{2} M}{t(T-t)}  \Big| \frac{\partial z}{\partial x} \Big|^{2}
+ \frac{C_{6} M^{3}}{(t(T-t))^{3}} |z|^{2} dxdt\, ,
\end{multline}
where $C_6=C_6(T,\beta):=C_4+C_3/2$. 
Coming back to our original variables thanks to identity (\ref{def-z}), we have
\begin{multline} \label{In-Carl-4}
\int_{0}^{T} \int\limits_{(-1,a') \cup (b',1)}  
\frac{C_{3} M^{3}|g|^{2} e^{-2\alpha}}{2(t(T-t))^{3}} dxdt
\\
\leqslant 
\int_{0}^{T} \int_{a'}^{b'} 
\left(
\frac{C_{7} M^{3} |g|^{2}}{(t(T-t))^{3}}
+  
\frac{C_{8} M}{t(T-t)}  
\Big| \frac{\partial g}{\partial x}\Big|^{2} 
\right)e^{-2\alpha} dxdt\, ,
\end{multline}
where $C_{8}=C_{8}(T,\beta):=2C_{2}$ and 
$C_{7}=C_{7}(T,\beta):=C_{6}+2C_{2} \sup \{ \beta'(x)^{2} ; x \in [a',b'] \}$.
Owing to (\ref{hyp-beta}) and the assumption $M \geqslant 1$, we have, for every $x \in [a,b]$, $t \in (0,T)$,
$$\frac{C_{7} M^{3}}{(t(T-t))^{3}} e^{-2\alpha}
\leqslant
\frac{C_{7} M^{3} }{(t(T-t))^{3}} e^{-\frac{2M}{t(T-t)}}
\leqslant
C_{9}\, ,$$
$$\frac{C_{8} M}{t(T-t)} e^{-2\alpha}
\leqslant
\frac{C_{8}t(T-t)}{M} \left( \frac{M}{t(T-t)} \right)^{2} e^{-\frac{2M}{t(T-t)}}
\leqslant
C_{10} t(T-t)\, ,$$
where
$C_{9}=C_{9}(T,\beta):=C_{7} \sup\{ x^{3}e^{-2x} ; x \in \mathbb{R}_{+} \}$
and
$C_{10}=C_{10}(T,\beta):=C_{8} \sup\{ x^{2} e^{-2x} ; x \in \mathbb{R}_{+} \}$.
Therefore, from (\ref{In-Carl-4}) we deduce
\begin{multline} \label{In-Carl-5}
\int_{0}^{T} \int\limits_{(-1,a') \cup (b',1)}  
\frac{C_{3} M^{3}|g|^{2} e^{-2\alpha}}{2(t(T-t))^{3}} dxdt
\\
\leqslant 
\int_{0}^{T} \int_{a'}^{b'} 
\left( C_{9} |g|^{2}  + 
C_{10} t(T-t) \Big| \frac{\partial g}{\partial x}\Big|^{2} \right) dxdt\,.
\end{multline}
Now, let us prove that the right hand side
of the previous inequality can be bounded by a first order term in $g$
on $(0,T) \times (a,b)$. 
We consider $\rho \in C^{\infty}(\mathbb{R},\mathbb{R}_{+})$ such that $0\le \rho \le 1$,
\begin{equation} \label{hyp-rho-hors(a,b)-carre}
\rho \equiv 1 \text{  on  } (a',b')\, ,
\end{equation}
\begin{equation} \label{hyp-rho-(a,b)-carre}
\rho \equiv 0 \text{  on  } (-1,a) \cup (b,1)\, .
\end{equation}
Multiplying the first equation of (\ref{adjointL2_n})
by $g \rho t(T-t)$ and then integrating over $(0,T) \times (-1,1 )$, we get
\begin{equation}\label{relation-rho-1-carre}
\int_{0}^{T}  \int_{-1}^1 \rho t(T-t)\left[
\frac{1}{2} \frac{d}{dt}\left( |g|^{2} \right)
- \frac{\partial^{2} g}{\partial x^{2}} g
+ (n\pi)^2 |x|^{2\gamma}|g|^2
\right] dx dt = 0\, .
\end{equation}
Integrating by parts with respect to space 
and time, we obtain
\begin{equation}\label{relation-rho-2-carre}
\frac{1}{2}\int_{0}^{T} \int_{-1}^1
 \frac{d}{dt} \Big[ |g|^{2} \Big] \rho t(T-t) 
dx dt
=
- \frac{1}{2}  \int_{0}^{T} \int_{-1}^1 |g|^{2} \rho (T-2t) 
dx dt\, ,
\end{equation}
\begin{multline} \label{relation-rho-3-carre}
-  \int_{0}^{T} \int_{-1}^1
\frac{\partial^{2} g}{\partial x^{2}} g \rho t(T-t) 
dx dt
\\ =
\int_{0}^{T}  \left[ \int_{-1}^1
\Big| \frac{\partial g}{\partial x} \Big|^{2} \rho t(T-t)
- \frac{1}{2} |g|^{2} \rho'' t(T-t) \right]
dx dt\, .
\end{multline}
Indeed, the boundary terms at $t=0$ and $t=T$ in (\ref{relation-rho-2-carre})
vanish owing to the factor $t(T-t)$, and the boundary terms at
$x=\pm 1$ in (\ref{relation-rho-3-carre})
vanish thanks to the boundary conditions on $g$.
Combining (\ref{relation-rho-1-carre}),
(\ref{relation-rho-2-carre}) and (\ref{relation-rho-3-carre}), we deduce
\begin{multline}\label{relation-rho-4-carre}
\int_{0}^{T} \int_{-1}^1  t(T-t)\left(
\Big| \frac{\partial g}{\partial x } \Big|^{2} \rho
- \frac{1}{2} |g|^{2} \rho''
+ (n\pi)^2 |x|^{2\gamma} |g|^2 \rho
\right) dx dt \\
- \frac{1}{2} \int_{0}^{T} \int_{-1}^1
|g|^{2} \rho (T-2t) dxdt = 0\, .
\end{multline}
In view of (\ref{hyp-rho-hors(a,b)-carre}), (\ref{hyp-rho-(a,b)-carre})
and (\ref{relation-rho-4-carre}), we have
\begin{multline}\label{relation-rho-4-carre*}
\int_{0}^{T} \int_{a'}^{b'}
\Big| \frac{\partial g}{\partial x } \Big|^{2} t(T-t) dx dt
\leqslant
\int_{0}^{T} \int_{-1}^1 
\Big| \frac{\partial g}{\partial x} \Big|^{2} \rho t(T-t) dx dt
\\ =
\frac{1}{2}\int_{0}^{T} \int_{-1}^1 
|g|^{2} [\rho (T-2t)  + t(T-t) (\rho'' - 2(n\pi)^2 |x|^{2\gamma} \rho) ]
dx dt
\\ \leqslant
C_{11} \int_{0}^{T} \int_{a}^b |g|^{2}  dx dt\, ,
\end{multline}
where 
$C_{11}=C_{11}(T,\rho):=T \|\rho\|_{L^{\infty}}
+ \frac{T^{2}}{2} \| \rho'' \|_{L^{\infty}}$.
Combining inequalities (\ref{relation-rho-4-carre*}) and (\ref{In-Carl-5}) leads to
\begin{equation} \label{In-Carl-6}
\int_{0}^{T} \int\limits_{(-1,a') \cup (b',1)}  
\frac{C_{3} M^{3}|g|^{2} e^{-2\alpha}}{(t(T-t))^{3}} dx dt
\leqslant 
\int_{0}^{T} \int_{a}^b
C_{12} |g|^{2} dx dt\, ,
\end{equation}
where 
$C_{12}=C_{12}(T,\beta,\rho):=2[C_{9}+C_{10}C_{11}]$.
Since
$$\frac{2T^{2}}{9} \leqslant t(T-t) \leqslant \frac{T^{2}}{4}\qquad \forall t \in 
\left[ \frac{T}{3}, \frac{2T}{3} \right]\,,$$
we have
$$\frac{e^{-2\alpha(t,x)}}{(t(T-t))^{3}}  
\geqslant
\frac{e^{-9c_{3}M/T^{2}}}{(T^{2}/4)^{3}} 
\quad
\forall x \in [-1,a'] \cup [b',1] \,,\  t \in (T/3,2T/3)\,,$$
where $c_{3}=c_3(\beta):=\sup\{ \beta(x) ; x \in [-1,a'] \cup [b',1] \}$. 
Therefore, (\ref{In-Carl-6}) implies
\begin{equation} 
\frac{C_{3}M^{3}}{(T^{2}/4)^{3}} e^{-\frac{9c_{3}M}{T^{2}}}
\int_{T/3}^{2T/3} \int\limits_{(-1,a') \cup (b',1)}  |g|^{2} dx dt
\leqslant 
\int_{0}^{T} \int_{a}^b
C_{12} |g|^{2} dx dt\, .
\end{equation}
Adding the same quantity to both sides
and using the inclusion $(a',b') \subset (a,b)$, we obtain
\begin{multline*} 
\frac{C_{3}M^{3}}{(T^{2}/4)^{3}} e^{-\frac{9c_{3}M}{T^{2}}}
\int_{T/3}^{2T/3} \int_{-1}^1  |g|^{2} dx dt \\
\leqslant 
\left( C_{12} + \frac{C_{3}M^{3}}{(T^{2}/4)^{3}} e^{-\frac{9c_{3}M}{T^{2}}} \right)
\int_{0}^{T} \int_{a}^b |g|^{2} dx dt\,,
\end{multline*}
which can also be written as
\begin{equation} 
\int_{T/3}^{2T/3} \int_{-1}^1  |g|^{2} dx dt
\leqslant 
\left( C_{12} \frac{(T^{2}/4)^{3}}{C_{3}M^{3}} e^{\frac{9c_{3}M}{T^{2}}} +1 \right)
\int_{0}^{T} \int_a^b |g|^{2} dx dt\, .
\end{equation}
Now, thanks to Proposition \ref{Prop:1st_eigenvalue},
$$ \int_{-1}^1 |g(T,x)|^2 dx 
\leqslant 
e^{-\frac{T}{3} c_* n^{\frac{2}{1+\gamma}}} \int_{-1}^{1} |g(t,x)|^2 dx\quad \forall t \in [T/3,2T/3]\, .$$
Thus,
\begin{equation} \label{g-dec}
\int_{-1}^1  |g(T,x)|^{2} dx
\leqslant 
\frac{3}{T} e^{-\frac{T}{3} c_* n^{\frac{2}{1+\gamma}}}
\left( C_{12} \frac{(T^{2}/4)^{3}}{C_{3}M^{3}} e^{\frac{9c_{3}M}{T^{2}}} +1 \right)
\int_{0}^{T} \int_a^b |g|^{2} dx dt\, .
\end{equation}
Now, let 
$$\mathcal{C}=\mathcal{C}(\beta):=\frac{1}{2} \sqrt{\frac{C_5}{2 C_3}}\,.$$
Then, there exists $n_1=n_1(T,\beta) \in \mathbb{N}^*$ such that,
for every $n \geqslant n_1$, the quantity $M=M(T,\beta,n)$ defined by (\ref{def:M}) satisfies
$$M=M(T,\beta,n)=\mathcal{C} n T^2\, .$$
Hence, for every $n \geqslant n_1$,
$$- c_* n^{\frac{2}{1+\gamma}} \frac{T}{3} + \frac{9c_{3}M}{T^{2}} = - c_* n^{\frac{2}{1+\gamma}} \frac{T}{3} + 9 c_3 \mathcal{C} n\, ,$$
where $c_*$, $c_3$, $\mathcal{C}$ depend only on $\beta$ and $\gamma$.

Let us assume that $\gamma \in [1/2,1)$ and $T>0$ is arbitrary. 
Since $2/(1+\gamma)>1$, there exists $n_2=n_2(\beta) \in \mathbb{N}^*$ such that
$$- c_* n^{\frac{2}{1+\gamma}} \frac{T}{3} + 9 c_3 \mathcal{C} n \leqslant 0\qquad \forall n \geqslant n_2\,.$$
So, inequality (\ref{g-dec}) yields, for every $n \geqslant \max\{n_1,n_2\}$,
\begin{equation} \label{final}
\int_{-1}^1  |g(T,x)|^{2} dx
\leqslant
\frac{3}{T}
\left( C_{12} \frac{(T^{2}/4)^{3}}{C_{3}M^{3}}  +1 \right)
\int_{0}^{T} \int_{a}^b |g|^{2} dx dt\, ,
\end{equation}
which in turn implies the conclusion.

Next, let us assume that $\gamma=1$ and $T>T_\sharp$, where
$$T_\sharp=T_\sharp(\beta):=\frac{27 c_3 \mathcal{C}}{c_*}\,.$$
Then, once again we recover (\ref{final}) for every $n \geqslant n_1$, and the conclusion follows as above.
\smallskip
\smallskip

\textbf{\underline{Second case: $\gamma \in (0,1/2)$.}}
\smallskip

The previous strategy does not apply to $\gamma \in (0,1/2)$ because the term $(n\pi)^2 [ |x|^{2\gamma} \alpha_x ]_x$ (that diverges at $x=0$) in (\ref{In-Carl-2}) can no longer be bounded by $\frac{C_{3} M^{3}}{(t(T-t))^{3}}$ (which is bounded at $x=0$).
Note that both terms are of the same order as $M^3$, because of the dependence of $M$ with respect to $n$ in (\ref{def:M}).
In order to deal with this difficulty, we adapt the choice of the weight $\beta$ and the dependence of $M$ with respect to $n$.
\smallskip

Let $\beta$ be a $C^1$-function on $(-1,1)$, which is also $C^2$ on $[-1,0)$ and $(0,1]$,
but such that $\beta''$ diverges at zero. More precisely, we assume that assumptions
(\ref{hyp-beta}), (\ref{hyp-beta'-(a,b)}) and (\ref{hyp_beta'-bord}) hold,
and condition (\ref{hyp-beta''}) is replaced by
\begin{equation}\label{hyp-beta''2}
\beta'' <0 \text{ on  }  [-1,0) \cup (0,a'] \cup [b',1]\,.
\end{equation}
Moreover, $\beta$ has the following form on a neighborhood $(-\epsilon,\epsilon)$ of $0$
\begin{equation} \label{hyp-beta-diverges}
\beta(x)=\mathcal{C}_0 - \int_0^x \sqrt{\text{sign}(s)|s|^{2\gamma}+\mathcal{C}_1}ds \qquad
\forall x \in (-\epsilon,\epsilon)\,,
\end{equation}
where the constant $\mathcal{C}_0$, $\mathcal{C}_1$ are large enough so that $\beta \geqslant 1$ and $\beta'<0$ on $(-\epsilon,\epsilon)$, respectively.
Then,
\begin{equation} \label{def:beta-gamma-pt}
\beta'(x)= - \sqrt{\text{sign}(x)|x|^{2\gamma}+\mathcal{C}_1}\qquad \forall x \in (-\epsilon,\epsilon)\, ,
\end{equation}
thus $\beta''$ diverges at $x=0$.
\smallskip

\noindent Performing the same computations as in the previous case, we get to inequality (\ref{In-Carl-1-BIS}).
Then, owing to (\ref{hyp-beta'-(a,b)}) and (\ref{hyp-beta''2}),
there exist $M_{1}=M_{1}(T,\beta)>0$, $C_{3}=C_{3}(\beta)>0$ and $C_{4}=C_{4}(T,\beta)>0$ such that, for every $M \geqslant M_{1}$ and $t \in (0,T)$,
\begin{multline*}
-\frac{1}{2}(\alpha_{t}-\alpha_{x}^{2})_{t} 
+  [ (\alpha_{t}-\alpha_{x}^{2}) \alpha_{x} ]_{x} 
- \frac{1}{2} \alpha_{xx}^{2}
\\
\geqslant
\frac{C_{3}M^{3}}{[t(T-t)]^{3}} |\beta''(x)|\, \beta'(x)^2\qquad
\forall x \in [-1,0) \cup (0,a'] \cup [b',1]\, , 
\end{multline*}
$$\Big| -\frac{1}{2}(\alpha_{t}-\alpha_{x}^{2})_{t} 
+  [ (\alpha_{t}-\alpha_{x}^{2}) \alpha_{x} ]_{x} 
- \frac{1}{2} \alpha_{xx}^{2} \Big|
\leqslant \frac{C_{4}M^{3}}{[t(T-t)]^{3}}
\qquad
\forall x \in [a',b']\, .$$
In view of \eqref{hyp_beta'-bord}, for every $M \geqslant M_{1}$,
\begin{multline}\label{In-Carl-2-gamma-pt} 
\int_{0}^{T} \int\limits_{(-1,a') \cup (b',1)}  
\left[
\frac{C_{3} M^{3}}{(t(T-t))^{3}} |\beta''(x)|\, \beta'(x)^2\,  |z|^{2}
+ (n\pi)^2 [ |x|^{2\gamma} \alpha_x ]_x |z|^2
\right] dxdt
\\
\leqslant 
\int_{0}^{T} \int_{a'}^{b'}
\left[
\frac{C_{2} M}{t(T-t)}  \Big| \frac{\partial z}{\partial x} \Big|^{2}
+ \frac{C_{4} M^{3}}{(t(T-t))^{3}} |z|^{2}
- (n\pi)^2 [ |x|^{2\gamma} \alpha_x ]_x |z|^2
\right] dxdt \, .
\end{multline} 
Moreover,
\begin{multline*}
|(n\pi)^2 [ |x|^{2\gamma} \alpha_x ]_x|  =
(n\pi)^2 \frac{M}{t(T-t)} \Big| 2\gamma \text{sign}(x) |x|^{2\gamma-1} \beta'(x) + |x|^{2\gamma} \beta''(x) \Big|
\\ \leqslant
\frac{C_5 n^2 M}{t(T-t)} \Big(  |x|^{2\gamma-1} |\beta'(x)| + |x|^{2\gamma} |\beta''(x)| \Big)\quad 
\forall x \in (-1,0) \cup (0,1)\, ,
\end{multline*}
where $C_5=\pi^2(2\gamma+1)$. From now on, we take
\begin{equation} \label{def:M-gamma-pt}
M=M(T,\beta,n):=\max \left\{ 1 ; M_1(T,\beta) ; \frac{n T^2}{\lambda} \right\}\, ,
\end{equation}
where $\lambda>0$ is a (small enough) constant, that will be chosen later on.

Then, there exists $n_1=n_1(T,\lambda,\beta) \in \mathbb{N}^*$ such that,
for every $n \geqslant n_1$, we have $M= nT^2/\lambda$. Therefore, for every $x \in (-1,0) \cup (0,1)$,
$$
|(n\pi)^2 [ |x|^{2\gamma} \alpha_x ]_x |
\leqslant
\frac{C_6 \lambda^2 M^3}{(t(T-t))^3} \Big(  |x|^{2\gamma-1} |\beta'(x)| + |x|^{2\gamma} |\beta''(x)| \Big)\, ,$$
where $C_6=C_6(\gamma)>0$. Let us verify that, for $\lambda>0$ small enough and for every $x \in (-1,0) \cup (0,a') \cup (b',1)$, we have
$$\begin{array}{c}
\displaystyle \frac{C_6 \lambda^2 M^3}{(t(T-t))^3}   |x|^{2\gamma-1} |\beta'(x)| 
\leqslant 
\frac{C_{3} M^{3}}{4(t(T-t))^{3}} |\beta''(x)|\, \beta'(x)^2\, ,
\\
\displaystyle \frac{C_6 \lambda^2 M^3}{(t(T-t))^3}  |x|^{2\gamma} |\beta''(x)| 
\leqslant 
\frac{C_{3} M^{3}}{4(t(T-t))^{3}} |\beta''(x)|\, \beta'(x)^2\, , 
\end{array}
$$
or, equivalently, for every $x \in (-1,0) \cup (0,a') \cup (b',1)$,
\begin{equation} \label{CN:beta-gamma-pt}
\begin{array}{c}
\displaystyle C_6 \lambda^2 |x|^{2\gamma-1}  \leqslant \frac{C_{3}}{4} |\beta''(x)|\cdot |\beta'(x)|\, , 
\\
\displaystyle C_6 \lambda^2  |x|^{2\gamma} \leqslant \frac{C_{3}}{4}  \beta'(x)^2\, .
\end{array}
\end{equation}
The second inequality is easy to satisfy (for $\lambda=\lambda(\beta,\gamma)$ small enough),
because $|\beta'|>0$ on $[-1,a'] \cup [b',1]$. Thanks to (\ref{def:beta-gamma-pt}), for every $x \in (-\epsilon,\epsilon)$,
$$\beta'(x)^2=\text{sign}(x)|x|^{2\gamma}+\mathcal{C}_1\, ,$$
so
$$\beta''(x) \beta'(x)= \gamma |x|^{2\gamma-1}\, .$$
Therefore, for every $x \in (-\epsilon,\epsilon)\setminus \{0\}$, the first inequality in (\ref{CN:beta-gamma-pt}) is equivalent to
$$C_6 \lambda^2   \leqslant \frac{C_{3}}{4} \gamma\, ,$$
which is trivially satisfied, when $\lambda=\lambda(\beta)$ is small enough.
Moreover, the first inequality of (\ref{CN:beta-gamma-pt}) holds for every $x \in [-1,-\epsilon]\cup[\epsilon,a']\cup[b',1]$
when $\lambda$ is small enough, since $|\beta''\beta'|>0$ on this compact set.
Finally, we deduce
\begin{multline} \label{In-Carl-3-gamma-pt}
\displaystyle \int_{0}^{T} \int\limits_{(-1,a') \cup (b',1)}  
\frac{C_{3} M^{3}}{2(t(T-t))^{3}} |\beta''(x)|\, \beta'(x)^2\,  |z|^{2} dxdt
\\
\displaystyle \leqslant 
\int_{0}^{T} \int_{a'}^{b'} \left[
\frac{C_{2} M}{t(T-t)}  \Big| \frac{\partial z}{\partial x} \Big|^{2}
+ \frac{C_{4} M^{3}}{(t(T-t))^{3}} |z|^{2} \right] dxdt\, ,
\end{multline} 
where $C_4=C_4(T,\beta)>0$. Since the function $|\beta''| (\beta')^2$ is bounded from below by some positive constant 
on $[-1,a']\cup [b',1]$, we also have
\begin{multline} \label{In-Carl-4-gamma-pt}
\int_{0}^{T} \int\limits_{(-1,a') \cup (b',1)}  
\frac{C_{3}' M^{3}}{2(t(T-t))^{3}}  |z|^{2} dxdt
\\
\leqslant 
\int_{0}^{T} \int_{a'}^{b'}
\left[
\frac{C_{2} M}{t(T-t)}  \Big| \frac{\partial z}{\partial x} \Big|^{2}
+ \frac{C_{4} M^{3}}{(t(T-t))^{3}} |z|^{2}
\right] dxdt\, ,
\end{multline} 
and the proof may be finished in the same way as in the first case. $\Box$

\subsection{Proof of the negative statements of Theorem \ref{thm:OI_n}}
\label{subsec:Proof_negative}

The goal of this section is the proof of the following results:
\begin{itemize}
\item if $\gamma=1$, then there exists $T_2>0$ such that, for every $T<T_2$, system (\ref{adjointL2_n}) is not observable in time $T$ uniformly with respect to $n$;
\item if $\gamma>1$ and $T>0$, then
system (\ref{adjointL2_n}) is not observable in time $T$ uniformly with respect to $n$.
\end{itemize}
The proof relies on the choice of particular test functions, that falsify uniform observability.
\smallskip

Let $\gamma \in [1,+\infty)$ be fixed and $T>0$.
For every $n \in \mathbb{N}^*$, we denote by  $\lambda_{n}$ (instead of $\lambda_{n,\gamma}$) the first
eigenvalue of the operator $A_{n,\gamma}$ defined in Section \ref{Dissipation speed},
and by $v_n$ the associated positive eigenvector of norm one, that is,
$$\left\lbrace \begin{array}{ll}
-v_n''(x)+[(n\pi)^2|x|^{2\gamma} -\lambda_n]v_n(x)=0\, ,\quad  x \in (-1,1)\,,\ n \in \mathbb{N}^*\,,\\
v_n(\pm 1)=0\,,\quad v_n\ge 0\, ,\\
\|v_n\|_{L^2(-1,1)}=1\,.
\end{array} \right.$$
Then, for every $n\ge 1$, the function
$$g_{n}(t,x):=v_n(x)e^{-\lambda_n t}\qquad \forall (t,x) \in \mathbb{R} \times (-1,1)\,,$$
solves the adjoint system (\ref{adjointL2_n}). Let us note that
$$\int_{-1}^1 g_n(T,x)^2 dx = e^{-2\lambda_n T}\, ,$$
$$\int_0^T \int_a^b  g_n(t,x)^2 dx dt = \frac{1-e^{-2\lambda_n T}}{2\lambda_n} \int_a^b v_n(x)^2 dx\, .$$
So, in order to prove that uniform observability fails, it suffices to show that
\begin{equation} \label{zero_cv}
\frac{e^{2\lambda_n T}}{\lambda_n}  \int_a^b v_n(x)^2 dx \rightarrow 0 \text{ when } n \rightarrow + \infty\, .
\end{equation}
In order to estimate the last integral, we will compare $v_n$ with an explicit supersolution of the problem on a suitable subinterval of $[-1,1]$.

\begin{Lem} \label{Lem:max_ppe}
Let $0<a<b<1$. For every $n \in \mathbb{N}^*$, set
\begin{equation} \label{def:xn}
x_n := \left( \frac{\lambda_n}{(n\pi)^2} \right)^{\frac{1}{2\gamma}}
\end{equation}
and let $W_n \in C^2([x_n,1],\mathbb{R})$ be a solution of
\begin{equation} \label{super-sol}
\left\lbrace \begin{array}{l}
-W_n''(x)+[(n\pi)^2 x^{2\gamma} -\lambda_n]W_n(x) \geqslant 0\, ,\quad  x \in (x_n,1)\,,\\
W_n(1) \geqslant 0\,,\\
W_n'(x_n) < -\sqrt{x_n} \lambda_n\,.
\end{array}\right.
\end{equation}
Then there exists $n_* \in \mathbb{N}^*$ such that, for every $n \geqslant n_*$,
$$\int_a^b v_n(x)^2 dx \leqslant \int_a^b W_n(x)^2 dx\,.$$ 
\end{Lem}

\noindent \textbf{Proof:} 
First, let us observe that, thanks to the second statement of Proposition~\ref{Prop:1st_eigenvalue}, $x_n \rightarrow 0$ when $n \rightarrow + \infty$. In particular, there exists $n_* \geqslant 1$ such that $x_n \leqslant a$ for every $n \geqslant n_*$. Now, let us prove that $|v_n'(x_n)| \leqslant \sqrt{x_n} \lambda_n$ for all \mbox{$n \geqslant n_*$.} Indeed, from Lemma \ref{Lem:firsteigeven}, we have $v_n(x)=v_n(-x)$, thus $v_n'(0)=0$. Hence, thanks to the Cauchy-Schwarz inequality and the relation $\|v_n\|_{L^2(-1,1)}=1$,
\begin{multline*}
\displaystyle |v_n'(x_n)| = \Big| \int_0^{x_n} v_n''(s) ds \Big| =
\Big| \int_0^{x_n} [(n\pi)^2|s|^{2\gamma}-\lambda_n]v_n(s) ds \Big|
\\ 
\displaystyle \leqslant
\left( \int_0^{x_n} [(n\pi)^2|s|^{2\gamma}-\lambda_n]^2 ds \right)^{1/2}
\left( \int_0^{x_n} v_n(s)^2 ds \right)^{1/2}
 \leqslant
\sqrt{x_n} \lambda_n\, .
\end{multline*}
Furthermore, we claim that $v_n(x) \leqslant W_n(x)$ for every $x \in [x_n,1]$, $n \geqslant n_*$. Indeed, if not so, there exists $x_* \in [x_n,1]$ such that
$$
(W_n-v_n)(x_*)=\min\{ (W_n-v_n)(x) ; x \in [x_n,1] \} <0\, .
$$
Since $(W_n-v_n)(1) \geqslant 0$ and $(W_n-v_n)'(x_n)<0$, we have $x_* \in (x_n,1)$.
Moreover, the function $W_n-v_n$ has a minimum at $x_*$, thus
$(W_n-v_n)'(x_*)=0$ and $(W_n-v_n)''(x_*) \geqslant 0$. Therefore,
$$-(W_n-v_n)''(x_*)+[(n\pi)^2 |x_*|^{2\gamma}-\lambda_n] (W_n-v_n)(x_*) <0\, ,$$
which is a contradiction. Our claim follows and the proof is complete. $\Box$
\smallskip

In order to apply Lemma \ref{Lem:max_ppe}, we look for an explicit 
supersolution $W_n$ of (\ref{super-sol}), of the form
\begin{equation} \label{def:W_n}
W_n(x)=C_n e^{-\mu_n x^{\gamma+1}}\, ,
\end{equation}
where $C_n, \mu_n>0$. Thus, the condition $W_n(1) \geqslant 0$ is automatically satisfied.
\smallskip

\emph{First step: Let us prove that, for an appropriate choice of $\mu_n$,
the first inequality of (\ref{super-sol}) holds.} Since
$$W_n'(x)=-\mu_n (\gamma+1) x^{\gamma} W_n(x)\,,$$
$$W_n''(x)= [ - \mu_n \gamma (\gamma+1) x^{\gamma-1} + \mu_n^2 (\gamma+1)^2 x^{2\gamma} ]W_n(x)\,,$$
the first inequality of (\ref{super-sol}) holds if and only if, for every $x \in (x_n,1)$,
\begin{equation} \label{rel1}
[ (n\pi)^2 - \mu_n^2 (\gamma+1)^2] x^{2\gamma} + \mu_n\gamma(\gamma+1)x^{\gamma-1} \geqslant \lambda_n\, .
\end{equation}
In particular, it holds when
\begin{equation} \label{mu_n:A1}
\mu_n \leqslant \frac{n\pi}{\gamma+1}
\end{equation}
and
\begin{equation} \label{mu_n:A2_1}
[ (n\pi)^2 - \mu_n^2 (\gamma+1)^2] x_n^{2\gamma} + \mu_n\gamma(\gamma+1)x_n^{\gamma-1} \geqslant \lambda_n\,.
\end{equation}
Indeed, in this case, the left hand side of (\ref{rel1}) is an increasing function of $x$. In view of (\ref{def:xn}), and after several simplifications, inequality (\ref{mu_n:A2_1}) can be recast as
$$\mu_n \leqslant \frac{\gamma}{\gamma+1} \left( \frac{(n\pi)^2}{\lambda_n} \right)^{\frac{1}{2} + \frac{1}{2\gamma}}\,.$$
So, recalling (\ref{mu_n:A1}), in order to satisfy the first inequality of (\ref{super-sol}) we can take
\begin{equation} \label{def:mu_n}
\mu_n := \min \left\{ 
\frac{n\pi}{\gamma+1} ; \frac{\gamma}{\gamma+1} \left( \frac{(n\pi)^2}{\lambda_n} \right)^{\frac{1}{2} + \frac{1}{2\gamma}}
\right\}\, .
\end{equation}
For the following computations, it is important to notice that, 
thanks to (\ref{def:mu_n}) and the second statement of Proposition \ref{Prop:1st_eigenvalue}, for $n$ large enough $\mu_n$ is of the form
\begin{equation} \label{mu_n:asymptotic}
\mu_n = C_1(\gamma) n\, .
\end{equation}

\emph{Second step: Let us prove that, for an appropriate choice of $C_n$, 
the third inequality of (\ref{super-sol}) holds.} Since
$$W_n'(x_n)=-C_n \mu_n (\gamma+1) x_n^\gamma e^{-\mu_n x_n^{\gamma+1}}\,,$$
the third inequality of (\ref{super-sol}) is equivalent to
$$C_n > \frac{\lambda_n e^{\mu_n x_n^{\gamma+1}}}{(\gamma+1) \mu_n x_n^{\gamma-\frac{1}{2}}}\,.$$
Therefore, it is sufficient to choose
\begin{equation} \label{def:C_n}
C_n:=\frac{2 \lambda_n e^{\mu_n x_n^{\gamma+1}}}{(\gamma+1) \mu_n x_n^{\gamma-\frac{1}{2}}}\, .
\end{equation}

\emph{Third step: Let us prove condition (\ref{zero_cv}).}
Thanks to Lemma \ref{Lem:max_ppe}, (\ref{def:W_n}), (\ref{mu_n:asymptotic}) and (\ref{def:C_n}), for every $n\geqslant n_*$,
\begin{multline*}
\frac{e^{2\lambda_n T}}{\lambda_n}  \int_a^b v_n(x)^2 dx
\leqslant
\frac{e^{2\lambda_n T}}{\lambda_n}  \int_a^b W_n(x)^2 dx
\leqslant
\frac{e^{2\lambda_n T}}{\lambda_n} W_n(a)^2
\\ \leqslant
\frac{e^{2\lambda_n T}}{\lambda_n} C_n^2 e^{-2\mu_n a^{1+\gamma}}
 \leqslant
\frac{e^{2\lambda_n T}}{\lambda_n}
\frac{4 \lambda_n^2 e^{2\mu_n x_n^{\gamma+1}}}{(\gamma+1)^2 \mu_n^2 x_n^{2\gamma-1}}
e^{-2\mu_n a^{1+\gamma}}.
\end{multline*}
By identities (\ref{def:xn}), (\ref{mu_n:asymptotic}) and Proposition \ref{Prop:1st_eigenvalue}, we have
$$\mu_n x_n^{\gamma+1} \leqslant C_2(\gamma)\qquad \forall n \in \mathbb{N}^*\,,$$
thus
\begin{equation} \label{last_bound}
\frac{e^{2\lambda_n T}}{\lambda_n}  \int_a^b v_n(x)^2 dx
\leqslant
e^{2n \left( \frac{\lambda_n}{n} T -  C_1(\gamma) a^{1+\gamma} \right)}
\frac{4 \lambda_n e^{2C_2(\gamma)}}{(\gamma+1)^2 \mu_n^2 x_n^{2\gamma-1}}\,.
\end{equation}

If $\gamma>1$, we deduce from the second statement of Proposition \ref{Prop:1st_eigenvalue} that
$$\frac{\lambda_n}{n} \rightarrow 0 \quad\text{ when } n \rightarrow + \infty\,.$$
So, for every $T>0$, there exists $n_\sharp \geqslant n_*$ such that, for every $n \geqslant n_\sharp$,
\begin{equation} \label{exp<0}
\frac{\lambda_n}{n} T -  C_1(\gamma) a^{1+\gamma} < - \frac{1}{2}C_1(\gamma) a^{1+\gamma}\,.
\end{equation}
Then, inequality (\ref{last_bound}) proves condition (\ref{zero_cv}) 
(since the term that multiplies the exponential behaves like a rational fraction of $n$).

If $\gamma=1$,  Proposition \ref{Prop:1st_eigenvalue} ensures that
$c_* n \leqslant \lambda_n \leqslant c^* n$, thus we deduce (\ref{exp<0}), hence (\ref{zero_cv}), for every
$$T<\frac{C_1(\gamma) a^{1+\gamma}}{2c_*}\, .$$

\subsection{End of the proof of Theorem \ref{thm:OI}} \label{subsec:ccl}

The first (resp. third) statement of Theorem \ref{thm:OI} has been proved 
in Subsection~\ref{subsec:Proof_positive} (resp. \ref{subsec:Proof_negative});
let us prove the second one.
\smallskip

Let us consider $\gamma=1$. Thanks to the results of Subsection \ref{subsec:Proof_positive}, the quantity
$$T^*:=\inf\{ T>0\, ;\ \text{system } (\ref{adjointL2_n}) \text{ is observable in time $T$ uniformly in } n \}$$
is well defined in $[0,+\infty)$. Moreover, as showed in Subsection \ref{subsec:Proof_negative}, $T^*>0$.
Clearly, uniform observability in some time $T_\sharp$ implies uniform observability in any time $T>T_\sharp$, so
\begin{itemize}
\item for every $T>T^*$, system (\ref{adjointL2_n}) is observable in time $T$ uniformly with respect to $n$;
\item for every $T<T^*$, system (\ref{adjointL2_n}) is not observable in time $T$ uniformly with respect to $n$.
\end{itemize}

\section{Conclusion and open problems}
\label{sec:ccl}

In this article we have studied the null controllability of the Grushin type equation (\ref{Grushin_eq}),
in a rectangle, with a distributed control localized on a strip parallel to the y-axis.
We have proved that null controllability
\begin{itemize}
\item holds in any positive time, when degeneracy is not too strong, i.e. $\gamma \in (0,1)$,
\item holds only in large time, when $\gamma=1$,
\item does not hold when degeneracy is too strong, i.e. $\gamma>1$.
\end{itemize}
Null controllability when $\gamma \in (0,1]$ and the control region $\omega$ is more general is an open problem.
When $\gamma=1$, it would be interesting to characterize the minimal time $T^*$ required 
for null controllability, and possibly connect it with the associated diffusion process.
Generalizations of this result to muldimensional configurations ($x \in (-1,1)^m$, $y \in (0,1)^n$), 
or boundary controls, are also open.

\bibliography{biblio3}
\bibliographystyle{plain}

\end{document}